\documentclass[12pt, reqno]{amsart}
\usepackage{eurosym}
\usepackage{amsmath, amsthm, amscd, amsfonts, amssymb, graphicx, color}
\usepackage[bookmarksnumbered, colorlinks, plainpages]{hyperref}

\hypersetup{colorlinks=true,linkcolor=red, anchorcolor=green, citecolor=cyan, urlcolor=red, filecolor=magenta, pdftoolbar=true}
\theoremstyle{definition}

\theoremstyle{remark}

\numberwithin{equation}{section}
\input{tcilatex}

\begin{document}
\title{New results on MS-Lipschitz summing operators}
\author{Maatougui Belaala, Athmane Ferradi and Khalil Saadi}
\address{Laboratory of Functional Analysis and Geometry of Spaces\\
University of M'sila\\
E-mail: \\
\ \ \ \ \ \ \ \ \ \ maatougui.belaala@univ-msila.dz\\
\ \ \ \ \ \ \ \ \ \ athmane.ferradi@univ-msila.dz \\
\ \ \ \ \ \ \ \ \ \ khalil.saadi@univ-msila.dz}
\date{}
\maketitle

\begin{abstract}
This paper focuses on the study of MS-Lipschitz $p$-summing operators, which
were initially defined by the authors in \cite{14}. Our objective is to
establish relationships between a Lipchitz mapping $T:X\rightarrow Y,$ where 
$X$ and $Y$ are pointed metric spaces, and its linearizations $\widehat{T},$ 
$\widetilde{T}$ and $T^{\#}$ for certain notion of summability. Moreover, we
extend our investigation by introducing new definitions in the category of
Lipschitz mappings defined on metric spaces, known as MS-strongly Lipschitz $%
p$-summing and MS-Lipschitz $p$-nuclear. We provide several results and
characterizations for these new concepts.%
%
%
%
%
%
\end{abstract}

\setcounter{page}{1}


\let\thefootnote\relax\footnote{
\par
{}\textit{Key words and phrases. }Strictly Lipschitz $p$-summing;
MS-Lipschitz $p$-summing; Lipschitz $p$-summing operators; MS-Lipschitz $p$%
-nuclear; $p$-summing operators; strongly $p$-summing operators; Pietsch
factorization; factorization theorems.}

\section{\textsc{Introduction and preliminaries}}

The theory of $p$-summing operators has undergone several stages of
development. It originated in the 1950s with Grothendieck's pioneering work 
\cite{8}, where he introduced the concept of $1$-summing operators. In 1967 
\cite{11}, Pietsch made a significant contribution by defining $p$-summing
operators for all positive values of $p$. His most notable result in this
theory is the Pietsch Factorization Theorem, which provides an integral
characterization. Since then, this theory has witnessed substantial
advancements. Researchers have expanded this theory in various directions,
including the sublinear, multilinear, and more recently, the Lipschitz case.
Let $Y$ be a pointed metric space and $G$ be a Banach space. We denote by $%
Lip_{0}\left( Y,G\right) $ the Banach space of all Lipschitz functions (or
Lipschitz operators) $T:Y\rightarrow G$ such that $T(0)=0$. This space is
equipped with pointwise addition and the Lipschitz norm. If $Y$ is equal to $%
\mathbb{K}$, we simply denoted by $Lip_{0}(Y)\left( =Y^{\#}\right) $. We
consider the evaluation functionals $\delta _{x}\in Lip_{0}(Y)^{\ast }$ for $%
x\in Y$ such that $\delta _{x}\left( s\right) =s\left( x\right) $ for $s\in
Lip_{0}(Y)$. The Lipschitz-free space $\mathcal{F}\left( Y\right) $ is the
closed space generated by these evaluation functionals. See \cite{6, 7, 10,
17, 18} for more details about free Banach spaces. We have $Y^{\#}=\mathcal{F%
}\left( Y\right) ^{\ast }$ holds isometrically. Every Lipschitz mapping $%
T:Y\rightarrow G$ induces a unique linear operator $\widehat{T}:\mathcal{F}%
\left( Y\right) \longrightarrow G$\ such that 
\begin{equation*}
\widehat{T}\circ \delta _{Y}=T,
\end{equation*}%
where $\delta _{Y}:Y\rightarrow \mathcal{F}\left( Y\right) $ is the
canonical embedding so that $\delta _{Y}\left( x\right) =\delta _{x}$\ for $%
x\in Y.$\ In this case, the identification 
\begin{equation}
Lip_{0}\left( Y,G\right) =\mathcal{B}\left( \mathcal{F}\left( Y\right)
,G\right) ,  \tag{1.1}
\end{equation}%
holds isometrically. It is well-known that if a Lipschitz mapping $%
T:Y\rightarrow G$ is Lipschitz $p$-summing, its linearization $\widehat{T}$
is not necessarily $p$-summing. However, in \cite[Theorem 3.5.]{14}, the
author introduced the concept of strictly Lipschitz $p$-summing, which
provides a well-established relation between $T$ and its linearization $%
\widehat{T}$ for the concept of $p$-summing. Now, let $Y$ and $W$ be two
pointed metric spaces. For every Lipschitz operator $T:Y\rightarrow W,$ we
can associate another linear operator $\widetilde{T}$\ such that the
following diagram commutes 
\begin{equation*}
\begin{array}{ccc}
Y & \underrightarrow{T} & W \\ 
\downarrow \delta _{Y} &  & \downarrow \delta _{W} \\ 
\mathcal{F}\left( Y\right)  & \underrightarrow{\widetilde{T}} & \mathcal{F}%
\left( W\right) 
\end{array}%
\end{equation*}%
i.e., $\widetilde{T}\circ \delta _{Y}=\delta _{W}\circ T$. The Lipschitz
adjoint map $T^{\#}:W^{\#}\rightarrow Y^{\#}$ of $T$ is defined as follows%
\begin{equation*}
T^{\#}\left( f\right) \left( x\right) =f\left( T\left( x\right) \right) ,%
\text{ for every }f\in W^{\#}\text{ and }x\in Y.
\end{equation*}%
Note that (as shown in \cite[p. 124]{7}), it is straightforward to verify
that%
\begin{equation}
\widetilde{T}^{\ast }=T^{\#},  \tag{1.2}
\end{equation}%
where $T^{\#}:W^{\#}\rightarrow Y^{\#}$ is the linear map called the
Lipschitz adjoint, which is defined by $T^{\#}\left( f\right) =f\circ T.$
Our objective in this paper is to establish relationships between a
Lipschitz mapping $T:Y\rightarrow W$, where $Y$ and $W$ are metric spaces,
and its linearizations $\widetilde{T}$ and $T^{\#}$. We focus on a specific
notion of summability and examine how the properties of $T$ are reflected in
its linearizations. In a previous work \cite{15}, the authors introduced the
concept of M-strictly Lipschitz $p$-summing operators, also referred to as
MS-Lipschitz $p$-summing, defined on pointed metric spaces and established a
significant connection between $T$ and its linearization $\widetilde{T}$.
Moreover, we will extend this idea by introducing a new definitions, called
MS-strongly Lipschitz $p$-summing\textit{\ }and\textit{\ }MS-Lipschitz $p$%
-nuclear, in the category of Lipschitz mappings defined on metric spaces.
The connection between $T$ and its linearizations, $\widetilde{T}$ and $%
T^{\#}$, will be established for these new concepts.

The structure of the paper is as follows:

Section 1 provides a brief review of the standard notations that will be
employed throughout the paper. Section 2 is dedicated to the examination of
various characterizations of MS-Lipschitz $p$-summing operators defined
between pointed metric spaces. In Section 3, we introduce the definition of
MS-strongly Lipschitz $p$-summing and MS-Lipschitz $p$-nuclear operators.
These operators exhibit remarkable properties, particularly their
associations with linearization operators. Furthermore, we present a Pietsch
Domination Theorem that applies specifically to these classes of operators.

Let's review some essential notations and terminology that will be used in
the following discussion. In this paper, we adopt the following conventions:
The letters $E$, $F$, and $G$ represent Banach spaces. The letters $Y$, $W$,
and $Z$ denote metric spaces with a distinguished point (\textit{pointed
metric spaces}), which we denote as 0. Given a metric space $Y$ and a Banach
space $G$, we represent their Lipschitz tensor product as $Y\boxtimes G$.
This tensor product is the vector space spanned by the linear functional $%
\delta _{\left( x,y\right) }\boxtimes h$ on $Lip_{0}\left( Y,G^{\ast
}\right) $, where $G^{\ast }$ denotes the topological dual of $G$. The
functional is defined as follows: For any $s\in Lip_{0}\left( Y,G^{\ast
}\right) $,%
\begin{equation*}
\delta _{\left( x,y\right) }\boxtimes h\left( s\right) =\left\langle s\left(
x\right) -s\left( y\right) ,h\right\rangle .
\end{equation*}%
When $\mathfrak{n}\in \mathcal{F}\left( Y\right) $ is expressed as $%
\mathfrak{n}=\sum_{j=1}^{m}\delta _{\left( x_{j},y_{j}\right) },$ we have
the following for every $s\in Y^{\#}$%
\begin{equation*}
\left\langle \mathfrak{n},s\right\rangle =\sum_{j=1}^{m}s(x_{j})-s(y_{j}).
\end{equation*}%
For a more detailed understanding of the properties of the space $Y\boxtimes
G$, we refer to \cite{1}. Now, let's consider a Lipschitz operator $%
T:Y\rightarrow W$ between pointed metric spaces, and $z=\sum_{l=1}^{k}\delta
_{\left( x_{l},y_{l}\right) }\boxtimes f_{l}\in Y\boxtimes W^{\#}.$ The
action of $T$ on $z$ is given by 
\begin{equation*}
\left\vert \left\langle T,z\right\rangle \right\vert =\left\vert
\sum_{l=1}^{k}f_{l}\left( T\left( x_{l}\right) \right) -f_{l}\left( T\left(
y_{l}\right) \right) \right\vert 
\end{equation*}%
Let $G$ be a Banach space. In the following definitions, $B_{G}$ represents
the closed unit ball of $G$, and $G^{\ast }$ denotes its (topological) dual.
For $1\leq p\leq \infty $ and $m\in \mathbb{N}^{\ast }$, we define two
Banach spaces as follows:

- $\ell _{p}^{m}\left( G\right) :$ This space consists of all sequences $%
\left( x_{j}\right) _{j=1}^{m}$ in $G$ with the norm given by%
\begin{equation*}
\left\Vert \left( x_{j}\right) _{j}\right\Vert _{\ell _{p}^{m}\left(
G\right) }=(\sum_{j=1}^{m}\left\Vert x_{j}\right\Vert ^{p})^{\frac{1}{p}}.
\end{equation*}

- $\ell _{p}^{m,w}\left( G\right) :$ This space comprises all sequences $%
\left( x_{j}\right) _{j=1}^{m}$ in $G,$ with the norm defined as%
\begin{equation*}
\left\Vert \left( x_{j}\right) _{j}\right\Vert _{\ell _{p}^{m,w}\left(
G\right) }=\underset{x^{\ast }\in B_{G^{\ast }}}{\sup }(\sum_{j=1}^{m}\left%
\vert \left\langle x_{j},x^{\ast }\right\rangle \right\vert ^{p})^{\frac{1}{p%
}}.
\end{equation*}%
If $G$ is equal to $\mathbb{K}$, we can simplify the notation as $\ell
_{p}^{m}$ and $\ell _{p}^{m,w}.$ In particular, let $\left( \mathfrak{n}%
_{j}\right) _{j=1}^{m_{1}}\in \mathcal{F}\left( Y\right) $ be such that $%
\mathfrak{n}_{j}=\sum_{i=1}^{m_{2}}\delta _{\left(
x_{j}^{i},y_{j}^{i}\right) }$ for $1\leq j\leq m_{1}.$\ In this case, we have%
\begin{eqnarray*}
\left\Vert \left( \mathfrak{n}_{j}\right) _{j}\right\Vert _{\ell
_{p}^{m,w}\left( \mathcal{F}\left( Y\right) \right) } &=&\sup_{s\in
Y^{\#}}(\sum_{j=1}^{m_{1}}\left\vert \left\langle \mathfrak{n}%
_{j},s\right\rangle \right\vert ^{p})^{\frac{1}{p}} \\
&=&\sup_{s\in Y^{\#}}(\sum_{j=1}^{m_{1}}\left\vert
\sum_{i=1}^{m_{2}}s(x_{j}^{i})-s(y_{j}^{i})\right\vert ^{p})^{\frac{1}{p}}.
\end{eqnarray*}%
Let us review the following concepts:

- The linear operator $u:F\rightarrow G$ is said to be $p$-summing if there
exists a positive constant $K$ such that, for any sequence $\left(
x_{j}\right) _{j=1}^{m}$ belonging to $F$, the following inequality holds%
\begin{equation}
\left\Vert \left( u\left( x_{j}\right) \right) _{j=1}^{m}\right\Vert _{\ell
_{p}^{m}\left( G\right) }\leq K\left\Vert \left( x_{j}\right)
_{j=1}^{m}\right\Vert _{\ell _{p}^{m,w}\left( F\right) }.  \tag{1.3}
\end{equation}%
The class of $p$-summing linear operators from the Banach space $F$ into $G$%
, denoted as $\Pi _{p}(F,G)$, forms a Banach space itself when equipped with
the norm $\pi _{p}(u)$. This norm is defined as the smallest constant $K$
for which the inequality (1.3) holds.\newline

- The linear operator $u:F\rightarrow G$ is said to be (Cohen) strongly $p$%
-summing if there exists a positive constant $K$ such that, for any sequence 
$\left( x_{j}\right) _{j=1}^{m}$ belonging to $F$ and any $\left(
y_{j}^{\ast }\right) _{j=1}^{m}$ belonging to $G^{\ast }$, the following
inequality holds%
\begin{equation}
\left\Vert \left( \left\langle u\left( x_{j},y_{j}^{\ast }\right)
\right\rangle \right) _{j=1}^{m}\right\Vert _{\ell _{1}^{m}}\leq K\left\Vert
\left( x_{j}\right) _{j=1}^{m}\right\Vert _{\ell _{p}^{m}\left( F\right)
}\left\Vert \left( y_{j}^{\ast }\right) _{j=1}^{m}\right\Vert _{\ell
_{p^{\ast }}^{m,w}\left( G^{\ast }\right) }.  \tag{1.4}
\end{equation}%
The class of strongly $p$-summing operators from the Banach space $F$ into $%
G $, denoted as $\mathcal{D}_{p}(F,G)$, forms a Banach space itself when
equipped with the norm $d_{p}(u)$. This norm is defined as the smallest
constant $K$ for which the inequality (1.4) holds.%
\vspace{0.5cm}%

- The linear operator $u:F\rightarrow G$ is said to be (Cohen) $p$-nuclear
if there exists a positive constant $K$ such that, for any sequence $\left(
x_{j}\right) _{j=1}^{m}$ belonging to $F$ and any $\left( y_{j}^{\ast
}\right) _{j=1}^{m}$ belonging to $G^{\ast }$, the following inequality holds%
\begin{equation}
\left\Vert \left( \left\langle u\left( x_{j},y_{j}^{\ast }\right)
\right\rangle \right) _{j=1}^{m}\right\Vert _{\ell _{1}^{m}}\leq K\left\Vert
\left( x_{j}\right) _{j=1}^{m}\right\Vert _{\ell _{p}^{m,w}\left( F\right)
}\left\Vert \left( y_{j}^{\ast }\right) _{j=1}^{m}\right\Vert _{\ell
_{p^{\ast }}^{m,w}\left( G^{\ast }\right) }.  \tag{1.5}
\end{equation}%
The class of strongly $p$-summing operators from the Banach space $F$ into $%
G $, denoted as $\mathcal{N}_{p}(F,G)$, forms a Banach space itself when
equipped with the norm $n_{p}(u)$. This norm is defined as the smallest
constant $K$ for which the inequality (1.5) holds.

\section{\textsc{M-strictly Lipschitz p-summing operators}}

Consider a pointed metric space $Y$ and let $G$ be a Banach space. The
concept of Lipschitz tensor product, denoted by $Y\boxtimes G$, was
introduced by Cabrera-Padilla et al. \cite{1}. An element $z$ in $Y\boxtimes
G$ can be represented as $z=\sum_{l=1}^{k}\delta _{\left( x_{l},y_{l}\right)
}\boxtimes h_{l}$ and can be viewed as a linear functional on $Lip_{0}\left(
Y,G^{\ast }\right) $. The action of this linear functional is defined by 
\begin{equation*}
\left\langle z,s\right\rangle =\sum_{l=1}^{k}\delta _{\left(
x_{l},y_{l}\right) }\boxtimes h_{l}\left( s\right) =\sum_{l=1}^{k}\left(
s\left( x_{l}\right) -s\left( y_{l}\right) \right) h_{l}\text{ for every }%
s\in Lip_{0}\left( Y,G^{\ast }\right) .
\end{equation*}%
The relationship between $Y\boxtimes G$ and $\mathcal{F}\left( Y\right)
\otimes G$ is straightforward, where $Y\boxtimes G$ is a vector subspace of $%
\mathcal{F}\left( Y\right) \otimes G$. Given an element $z\in Y\boxtimes G$,
we can define the set $A_{z}$ as the set of all representations of $z$ in $%
\mathcal{F}\left( Y\right) \otimes G$, that is, 
\begin{equation}
A_{z}:=\left\{ \left( \left( \mathfrak{n}_{j}\right) _{j=1}^{m},\left(
h_{j}\right) _{j=1}^{m}\right) :m\in \mathbb{N},\text{ }\mathfrak{n}_{j}\in 
\mathcal{F}\left( Y\right) ,h_{j}\in G:z=\sum_{j=1}^{m}\mathfrak{n}%
_{j}\otimes h_{j}\right\} .  \tag{2.1}
\end{equation}

Let's consider $\beta $ as a tensor norm defined on Banach spaces. Based on 
\cite[Theorem 3.1]{14}, it has been established that there exists a
corresponding Lipschitz cross-norm, denoted as $\beta ^{L},$ defined on $%
Y\boxtimes G$ by: 
\begin{equation}
\beta ^{L}(\sum_{l=1}^{k}\delta _{\left( x_{l},y_{l}\right) }\boxtimes
h_{l})=\beta (\sum_{l=1}^{k}\delta _{\left( x_{l},y_{l}\right) }\otimes
h_{l}),  \tag{2.2}
\end{equation}%
with $\sum_{l=1}^{k}\delta _{\left( x_{l},y_{l}\right) }\otimes h_{l}\ $is
an element of $\mathcal{F}\left( Y\right) \otimes G$. Before presenting the
following definition, it is necessary to recall the norms of Chevet-Saphar $%
d_{p}$ and $g_{p}$ \cite{2, 12, 16}, which are defined on two Banach spaces $%
F$ and $G$%
\begin{equation*}
d_{p}\left( z\right) =\inf \left\{ \left\Vert \left( x_{j}\right)
_{j=1}^{m}\right\Vert _{\ell _{p^{\ast }}^{m,w}\left( F\right) }\left\Vert
\left( y_{j}\right) _{j=1}^{m}\right\Vert _{\ell _{p}^{m}\left( G\right)
}\right\} \text{ },
\end{equation*}%
and 
\begin{equation*}
g_{p}\left( z\right) =\inf \left\{ \left\Vert \left( x_{j}\right)
_{j=1}^{m}\right\Vert _{\ell _{p^{\ast }}^{m}\left( F\right) }\left\Vert
\left( y_{j}\right) _{j=1}^{m}\right\Vert _{\ell _{p}^{m,w}\left( G\right)
}\right\} ,
\end{equation*}%
here the $\inf $ is taken over all representations of $z$ in the form of $%
z=\sum_{j=1}^{m}x_{j}\otimes y_{j}\in F\otimes G.$ It is worth noting that
we can utilize the Chevet-Saphar norms to provide equivalent definitions for
(1.3) and (1.4) (see \cite[p. 140]{12}). Specifically, the linear operator $%
u:F\rightarrow G$ is said to be $p$-summing if there exists a constant $K>0$
such that for any $z=\sum_{j=1}^{m}x_{j}\otimes y_{j}^{\ast }\in F\otimes
G^{\ast },$ the following inequality holds:%
\begin{equation*}
\left\vert \left\langle u,z\right\rangle \right\vert =\left\vert
\sum_{j=1}^{m}\left\langle u\left( x_{j}\right) ,y_{j}^{\ast }\right\rangle
\right\vert \leq Kd_{p}\left( z\right) .
\end{equation*}%
If we replace $d_{p}\left( u\right) $ with $g_{p}\left( u\right) $, we
obtain the definition of strongly $p$-summing. In \cite{14}, the Lipschitz
cross-norm $d_{p}^{L}$ is defined as follows:%
\begin{equation*}
d_{p}^{L}\left( z\right) =d_{p}(\sum_{l=1}^{k}\delta _{\left(
x_{l},y_{l}\right) }\otimes t_{l})=\inf_{\left( \left( \mathfrak{n}%
_{j}\right) _{j=1}^{m},\left( h_{j}\right) _{j=1}^{m}\right) \in
A_{z}}\left\{ \left\Vert \left( \mathfrak{n}_{j}\right)
_{j=1}^{m}\right\Vert _{l_{p}^{m,w}\left( \mathcal{F}\left( Y\right) \right)
}\left\Vert \left( h_{j}\right) _{j=1}^{m}\right\Vert _{l_{p^{\ast
}}^{m}\left( G\right) }\right\} .
\end{equation*}

Most definitions of summability for Lipschitz mappings are typically defined
from a metric space to a Banach space. However, in the following definition
of MS-Lipschitz $p$-summing, we consider Lipschitz operators defined on
metric spaces. This new perspective allows us to establish a meaningful
relationship between $T$ and its linearization $\widetilde{T}$. Before
proceeding, let us recall the following definition introduced in \cite{14}.%
\newline

\textbf{Definition 2.1. }\cite{14} Consider $1\leq p\leq \infty ,$ $Y$ be a
pointed\textit{\ }metric space and $G$ be a Banach space. The Lipschitz
operator $T:Y\rightarrow G$ is considered to be\textit{\ strictly Lipschitz }%
$p$\textit{-summing} if there exists a constant $K>0$ such that for every $%
z=\sum_{l=1}^{k}\delta _{\left( x_{l},y_{l}\right) }\boxtimes h_{l}^{\ast
}\in Y\boxtimes G^{\ast }$ we have%
\begin{equation}
\left\vert \left\langle T,z\right\rangle \right\vert \leq Kd_{p}^{L}(z). 
\tag{2.3}
\end{equation}

Building upon the aforementioned idea, we have introduced the concept of
MS-Lipschitz $p$-summing operators \cite{15}. In contrast to considering
elements from the dual space of $G$, we now focus on elements from its
Lipschitz space $G^{\#}$.\newline

\textbf{Definition 2.2. }\cite{15} For $1\leq p\leq \infty $ and $Y,W$ being
two pointed\textit{\ }metric spaces, a Lipschitz operator $T:Y\rightarrow W$
is considered to be MS-Lipschitz $p$-summing if there exists a constant $K>0$
such that for every $z=\sum_{l=1}^{k}\delta _{\left( x_{l},y_{l}\right)
}\boxtimes f_{l}\in Y\boxtimes W^{\#}$ we have%
\begin{equation}
\left\vert \left\langle T,z\right\rangle \right\vert \leq Kd_{p}^{L}(z). 
\tag{2.4}
\end{equation}

The set of all MS-Lipschitz $p$-summing operators from $Y$ into $W$ is
denoted by $\Pi _{p}^{MSL}\left( Y,W\right) $. It represents the collection
of operators that satisfy the MS-Lipschitz $p$-summing property. The
constant $\pi _{p}^{MSL}\left( T\right) $ corresponds to the smallest value
of $K$ that satisfies inequality (2.4) for a given operator $T$. It is
important to note that if $W$ is a Banach space, the set $\Pi
_{p}^{MSL}\left( Y,W\right) $ does not possess the structure of a vector
space.\newline

\textbf{Proposition 2.3}. \textit{Let }$1\leq p\leq \infty .$\textit{Every
MS-Lipschitz }$p$\textit{-summing operator from a pointed metric space }$Y$%
\textit{\ into a Banach space }$G$\textit{\ is strictly Lipschitz }$p$%
\textit{-summing.}\newline

\textbf{Proof}. Let $T:X\rightarrow E$ be a MS-Lipschitz $p$-summing
operator. Let $x_{l},y_{l}\in Y$ and $h_{k}^{\ast }\in E^{\ast }$ $\left(
1\leq l\leq k\right) ,$ then%
\begin{eqnarray*}
\left\vert \left\langle T,z\right\rangle \right\vert  &=&\left\vert
\sum_{l=1}^{k}h_{l}^{\ast }\left( T\left( x_{l}\right) \right) -h_{l}^{\ast
}\left( T\left( y_{l}\right) \right) \right\vert  \\
&\leq &Kd_{p}^{L}(z),
\end{eqnarray*}%
where $z=\sum_{l=1}^{k}\delta _{\left( x_{l},y_{l}\right) }\boxtimes
h_{k}^{\ast }\in Y\boxtimes G^{\#}.$ So, as $Y\boxtimes G^{\ast }\subset
Y\boxtimes G^{\#},$ the definition of $d_{p}^{L}(z)$ on $Y\boxtimes G^{\#}$
is smaller than on $Y\boxtimes G^{\ast },$ consequently, the condition $%
\left( 2.3\right) $ is verified.$\qquad \blacksquare $\newline

\textbf{Remark 2.4}. In the case where $F$ and $G$ are Banach spaces, it is
well-known that the definitions of strictly Lipschitz $p$-summing, Lipschitz 
$p$-summing, and $p$-summing coincide for linear operators from $F$ to $G$
(see \cite[Proposition 3.8]{14}). Furthermore, in our specific case, the
definition of MS-Lipschitz $p$-summing implies $p$-summing; however, the
converse is not true, as illustrated in the following example: Consider the
identity operator $id_{F}:F\rightarrow F$. It can be easily demonstrated
that $\widetilde{id_{F}}=id_{\mathcal{F}\left( F\right) }$, indicating that
the following diagram is commutative%
\begin{equation*}
\begin{array}{ccc}
F & \underrightarrow{id_{F}} & F \\ 
\downarrow \delta _{F} &  & \downarrow \delta _{F} \\ 
\mathcal{F}\left( F\right)  & \underrightarrow{id_{\mathcal{F}\left(
F\right) }} & \mathcal{F}\left( F\right) 
\end{array}%
\end{equation*}%
If $F$ is a finite-dimensional space, then $id_{F}$ is indeed $p$-summing,
and consequently, it is also strictly Lipschitz $p$-summing. However, $id_{%
\mathcal{F}\left( F\right) }$ cannot be $p$-summing since $\mathcal{F}\left(
F\right) $ is not finite-dimensional. Therefore, $id_{F}$ is not
MS-Lipschitz $p$-summing.\newline

The following statement presents the main result of this section.\newline

\textbf{Theorem 2.5}.\textit{\ Consider }$1\leq p\leq \infty .$\textit{\ Let 
}$Y$\textit{\ and }$W$\textit{\ be two pointed metric spaces.} \textit{Let }$%
T:X\rightarrow Y$\textit{\ be a Lipschitz operator. The following properties
are equivalent.}

\noindent \textit{1) }$T$\textit{\ is MS-Lipschitz }$p$\textit{-summing.}

\noindent \textit{2)} $\widetilde{T}:\mathcal{F}\left( Y\right) \rightarrow 
\mathcal{F}\left( W\right) $\textit{\ is }$p$\textit{-summing.}

\noindent \textit{3)} $\delta _{Y}\circ T:Y\rightarrow \mathcal{F}\left(
W\right) $ \textit{is strictly Lipschitz }$p$\textit{-summing.}

\noindent \textit{4) There is a constant }$K>0$\textit{\ such that for every 
}$\left( x_{j}^{i}\right) _{j=1}^{m_{1}},\left( y_{j}^{i}\right)
_{j=1}^{m_{1}}$\textit{\ in }$Y;\left( 1\leq i\leq m_{2}\right) $ \textit{and%
} $m_{1},m_{2}\in \mathbb{N}^{\ast }$\textit{, we have}%
\begin{equation}
(\sum_{j=1}^{m_{1}}\left\Vert \sum_{i=1}^{m_{2}}\delta _{W}\circ
T(x_{j}^{i})-\delta _{W}\circ T(y_{j}^{i})\right\Vert ^{p})^{\frac{1}{p}%
}\leq K\sup_{s\in Y^{\#}}(\sum_{j=1}^{m_{1}}\left\vert
\sum_{i=1}^{m_{2}}s(x_{j}^{i})-s(y_{j}^{i})\right\vert ^{p})^{\frac{1}{p}}. 
\tag{2.5}
\end{equation}

\textbf{Proof}. $1)\Leftrightarrow 2)$ See \cite[Proposition 2.4.]{15}$.$

\noindent $2)\Rightarrow 3):$ Suppose that $\widetilde{T}$\textit{\ }is $p$%
-summing. Then%
\begin{eqnarray*}
\sum_{j=1}^{m_{1}}(\left\Vert \sum_{i=1}^{m_{2}}\delta _{W}\circ
T(x_{j}^{i})-\delta _{W}\circ T(y_{j}^{i})\right\Vert ^{p})^{\frac{1}{p}}
&=&(\sum_{j=1}^{m_{1}}\left\Vert \widetilde{T}\left( \mathfrak{n}_{j}\right)
\right\Vert ^{p})^{\frac{1}{p}} \\
&\leq &\pi _{p}\left( \widetilde{T}\right) \sup_{s\in
Y^{\#}}(\sum_{j=1}^{m_{1}}\left\vert s\left( \mathfrak{n}_{j}\right)
\right\vert ^{p})^{\frac{1}{p}} \\
&\leq &\pi _{p}\left( \widetilde{T}\right) \sup_{s\in
Y^{\#}}(\sum_{j=1}^{m_{1}}\left\vert
\sum_{i=1}^{m_{2}}s(x_{j}^{i})-s(y_{j}^{i})\right\vert ^{p})^{\frac{1}{p}}
\end{eqnarray*}%
$3)\Rightarrow 2):$ We have%
\begin{eqnarray*}
(\sum_{j=1}^{m_{1}}\left\Vert \widetilde{T}\left( \mathfrak{n}_{j}\right)
\right\Vert ^{p})^{\frac{1}{p}} &=&(\sum_{j=1}^{m_{1}}\left\Vert
\sum_{i=1}^{m_{2}}\delta _{W}\circ T\left( x_{j}^{i}\right) -\delta
_{W}\circ T\left( y_{j}^{i}\right) \right\Vert ^{p})^{\frac{1}{p}} \\
&\leq &K\sup_{s\in Y^{\#}}(\sum_{j=1}^{m_{1}}\left\vert
\sum_{i=1}^{m_{2}}s(x_{j}^{i})-s(y_{j}^{i})\right\vert ^{p})^{\frac{1}{p}} \\
&\leq &K\sup_{s\in Y^{\#}}(\sum_{j=1}^{m_{1}}\left\vert s\left( \mathfrak{n}%
_{j}\right) \right\vert ^{p})^{\frac{1}{p}}.
\end{eqnarray*}%
As a consequence, $\widetilde{T}$ is $p$-summing, and by applying the result
in \cite[Theorem 2.12]{4}, we can obtain the desired result.

\noindent $3)\Leftrightarrow 4)$ It is immediate.$\qquad \blacksquare $%
\newline

By setting $m_{2}=1$ in formula (2.5) and considering the isometric property
of $\delta _{W}$, we arrive at the precise formulation of Lipschitz $p$%
-summing as originally defined by Farmer \cite{5}, indeed 
\begin{eqnarray*}
(\sum_{j=1}^{m_{1}}\left\Vert \delta _{W}\circ T(x_{j})-\delta _{W}\circ
T(y_{j})\right\Vert ^{p})^{\frac{1}{p}} &=&(\sum_{j=1}^{m_{1}}d\left(
T(x_{j}),T(y_{j})\right) ^{p})^{\frac{1}{p}} \\
&\leq &K\sup_{s\in Y^{\#}}(\sum_{j=1}^{m_{1}}\left\vert
s(x_{j})-s(y_{j})\right\vert ^{p})^{\frac{1}{p}}.
\end{eqnarray*}

\textbf{Corollary 2.6}. \textit{Consider a Lipschitz operator }$%
T:Y\rightarrow W$\textit{\ between pointed metric spaces. The following
properties are equivalent.}

\noindent \textit{1) }$T$\textit{\ is MS-Lipschitz }$p$\textit{-summing.}

\noindent \textit{2) The Lipschitz adjoint }$T^{\#}:Y^{\#}\rightarrow X^{\#}$
\textit{is strongly }$p^{\ast }$\textit{-summing.}\newline

\textbf{Proof}. According to (1.2), the dual operator of $\widetilde{T}$ is $%
T^{\#}$. To establish the equivalence between the given properties, we can
make use of the result mentioned in \cite[Theorem 2.2.2]{3}.\qquad $%
\blacksquare $\newline

\textbf{Proposition 2.7}. \textit{Consider a Lipschitz operator }$%
T:Y\rightarrow W$\textit{\ between pointed metric spaces such that }$Y$%
\textit{\ or }$W$\textit{\ is finite, then }$T$\textit{\ is MS-Lipschitz }$p$%
\textit{-summing.}\newline

\textbf{Proof}. Suppose $Y$ is a finite metric space. According to \cite[%
Example 2.3.6]{17}, we know that the space $\mathcal{F}\left( Y\right) $ is
finite-dimensional. Therefore, the linearization $\widetilde{T}:\mathcal{F}%
\left( Y\right) \rightarrow \mathcal{F}\left( W\right) $\textit{\ }is $p$%
-summing, and consequently, $T:Y\rightarrow W$ is MS-Lipschitz $p$%
-summing.\qquad $\blacksquare $\newline

The Pietsch domination theorem is an intriguing characterization that is
satisfied by the class of MS-Lipschitz $p$-summing operators.\newline

\textbf{Theorem 2.8}. \textit{Consider a Lipschitz operator }$T:Y\rightarrow
W$\textit{\ between pointed metric spaces. We have the following equivalent
properties.}

\noindent \textit{1) }$T$\textit{\ is MS-Lipschitz }$p$\textit{-summing.}

\noindent \textit{2) There exist a constant }$K>0$, \textit{a} \textit{Radon}
\textit{probability }$\mu $\textit{\ on }$B_{Y^{\#}}$ \textit{such that for
every }$\left( x^{i}\right) _{i=1}^{m},\left( y^{i}\right) _{i=1}^{m}\subset
Y,$ \textit{we have}%
\begin{equation}
\left\Vert \sum_{i=1}^{m}\delta _{W}\circ T(x^{i})-\delta _{W}\circ
T(y^{i})\right\Vert \leq K(\int_{B_{Y^{\#}}}\left\vert \sum_{i=1}^{m}s\left(
x^{i}\right) -s\left( y^{i}\right) \right\vert ^{p}d\mu \left( s\right) )^{%
\frac{1}{p}}.  \tag{2.6}
\end{equation}%
\textit{In this case, we have}%
\begin{equation*}
\pi _{p}^{MSL}\left( T\right) =\inf \left\{ K:\text{ \textit{verifying} 
\textit{(2.6)}}\right\} .
\end{equation*}

\textbf{Proof}. $1)\Rightarrow 2):$ Since $T$\textit{\ }is MS-Lipschitz $p$%
-summing,\textit{\ }then $\widetilde{T}:\mathcal{F}\left( Y\right)
\rightarrow \mathcal{F}\left( W\right) $ is $p$-summing. By Pietsch
domination theorem for $p$-summing linear operator \cite[Theorem 2.12]{4},
we have%
\begin{eqnarray*}
\left\Vert \widetilde{T}\left( \sum_{i=1}^{m}\delta _{\left(
x^{i},y^{i}\right) }\right) \right\Vert  &\leq &\pi _{p}\left( \widetilde{T}%
\right) (\int_{B_{Y^{\#}}}\left\vert \left\langle \sum_{i=1}^{m}\delta
_{\left( x^{i},y^{i}\right) },s\right\rangle \right\vert ^{p}d\mu \left(
s\right) )^{\frac{1}{p}} \\
&\leq &\pi _{p}\left( \widetilde{T}\right) (\int_{B_{Y^{\#}}}\left\vert
\sum_{i=1}^{m}s\left( x^{i}\right) -s\left( y^{i}\right) \right\vert
^{p}d\mu \left( s\right) )^{\frac{1}{p}}.
\end{eqnarray*}%
On the other hand, 
\begin{eqnarray*}
\left\Vert \widetilde{T}\left( \sum_{i=1}^{m}\delta _{\left(
x^{i},y^{i}\right) }\right) \right\Vert  &=&\left\Vert \sum_{i=1}^{m}%
\widetilde{T}\left( \delta _{\left( x^{i},y^{i}\right) }\right) \right\Vert 
\\
&=&\left\Vert \sum_{i=1}^{m}\widetilde{T}\circ \delta _{Y}\left(
x^{i}\right) -\widetilde{T}\circ \delta _{Y}\left( y^{i}\right) \right\Vert 
\\
&=&\left\Vert \sum_{j=1}^{n}\delta _{W}\circ T(x^{i})-\delta _{W}\circ
T(y^{i})\right\Vert 
\end{eqnarray*}%
Therefore, we have obtained the desired result.

$2)\Rightarrow 1):$ Similarly, we can apply the same argument.\qquad $%
\blacksquare $\newline

We conclude this section with a result concerning Lipschitz operators that
have a finite image. The following Lemma establishes a relationship between
the free space of $T\left( Y\right) $ and the image $\widetilde{T}\left( 
\mathcal{F}\left( Y\right) \right) $.\newline

\textbf{Lemma 2.9}. \textit{Let }$X$\textit{\ and }$Y$\textit{\ be two
pointed metric spaces. Consider a Lipschitz operator }$T:Y\rightarrow W$ 
\textit{such that} $T\left( Y\right) $\textit{\ is a closed subset of }$W$%
\textit{. Then, we have the following}%
\begin{equation*}
\widetilde{T}\left( \mathcal{F}\left( Y\right) \right) =\mathcal{F}\left(
T\left( Y\right) \right) .
\end{equation*}

\textbf{Proof}. By \cite[Theorem 2.2.6]{17}, we have%
\begin{eqnarray*}
\mathcal{F}\left( T\left( Y\right) \right)  &=&\overline{span}\left\{ \delta
_{T\left( x\right) }:x\in Y\right\}  \\
&=&\overline{span}\left\{ \delta _{W}\left( T\left( x\right) \right) :x\in
Y\right\}  \\
&=&\overline{span}\left\{ \widetilde{T}\left( \delta _{x}\right) :x\in
Y\right\}  \\
&=&\widetilde{T}\left( \mathcal{F}\left( Y\right) \right) .\qquad
\blacksquare 
\end{eqnarray*}

\textbf{Corollary 2.10}.\textit{\ Let }$Y$\textit{\ and }$W$\textit{\ be two
pointed metric spaces. Suppose that }$T:Y\rightarrow W$\textit{\ is a
Lipschitz operator such that }$T\left( Y\right) $\textit{\ is a finite set.
Then, the linearization }$\widetilde{T}$\textit{\ has finite rank.
Consequently, we can conclude that every finite Lipschitz operator is
MS-Lipschitz }$p$\textit{-summing}.

\section{\textsc{Cohen MS-Lipschitz }$p$\textsc{-nuclear operators}}

In \cite{3}, Cohen introduced the concepts of strongly $p$-summing and $p$%
-nuclear operators in the category of linear operators. Since then, many
authors have explored and extended these notions in various directions,
including multilinear, sublinear, and Lipschitz cases. We will further
extend these concepts using a similar approach to the one presented in the
previous section.\newline

\textbf{Definition 3.1}. Consider a Lipschitz operator $T:Y\rightarrow W$\
between pointed metric spaces. For $1\leq p\leq \infty ,$ $T$ is said to be%
\textit{\ MS-strongly Lipschitz }$p$\textit{-summing} if there exists a
constant $K>0$ such that the following condition holds for any $%
z=\sum_{l=1}^{k}\delta _{\left( x_{l},y_{l}\right) }\boxtimes f_{l}\in
Y\boxtimes W^{\#}$%
\begin{equation}
\left\vert \left\langle T,z\right\rangle \right\vert \leq Kg_{p}^{L}(z). 
\tag{3.1}
\end{equation}%
We denote the set of all MS-strongly\textit{\ }Lipschitz $p$-summing
operators from $Y$ into $W$ as $\mathcal{D}_{p}^{MSL}\left( Y,W\right) ,$
and $d_{p}^{MSL}\left( T\right) $ represents the smallest constant $K$ that
satisfies (3.1). If $W$ is a Banach space, it's important to note that $%
\mathcal{D}_{p}^{MSL}\left( Y,W\right) $ does not possess the structure of a
vector space.\newline

\textbf{Theorem 3.2.} \textit{Consider }$1\leq p\leq \infty $\textit{\ and a
Lipschitz operator }$T:Y\rightarrow W$\textit{\ between pointed metric
spaces. The following properties are equivalent.}

\noindent \textit{1) }$T$\textit{\ is MS-strongly Lipschitz }$p$\textit{%
-summing.}

\noindent \textit{2) There exists a constant }$K>0$\textit{\ such that} 
\textit{for every }$\left( x_{j}\right) _{j=1}^{m},\left( y_{j}\right)
_{j=1}^{m}$\textit{\ in }$Y$ \textit{and} $\left( f_{j}\right) _{j=1}^{m}$%
\textit{in} $W^{\#};$ \textit{(}$m\in \mathbb{N}^{\ast }$\textit{), }%
\begin{equation}
\left\vert \sum_{j=1}^{m}f_{j}\left( T(x_{j})\right) -f_{j}\left(
T(y_{j})\right) \right\vert \leq K(\sum_{j=1}^{m}d\left( x_{j},y_{j}\right)
^{p})^{\frac{1}{p}}\left\Vert \left( f_{j}\right) _{j=1}^{m}\right\Vert
_{\ell _{p^{\ast }}^{m,w}\left( W^{\#}\right) }.  \tag{3.2}
\end{equation}

\noindent \textit{3) }$\delta _{W}\circ T:Y\rightarrow \mathcal{F}\left(
W\right) $ \textit{is strongly Lipschitz }$p$\textit{-summing.}

\noindent \textit{4)} $\widetilde{T}:\mathcal{F}\left( Y\right) \rightarrow 
\mathcal{F}\left( W\right) $\textit{\ is strongly }$p$\textit{-summing.}%
\newline

\textbf{Proof}.

\noindent $1)\Rightarrow 2):$ Let $T$ be a MS-strongly Lipschitz $p$-summing
operator. Let\textit{\ }$\left( x_{j}\right) _{j=1}^{m},\left( y_{j}\right)
_{j=1}^{m}$\textit{\ }in\textit{\ }$Y$ and $\left( f_{j}\right)
_{j=1}^{m}\subset W^{\#}.$ We get%
\begin{equation*}
\left\vert \left\langle T,z\right\rangle \right\vert =\left\vert
\sum_{j=1}^{m}f_{j}\left( T(x_{j})\right) -f_{j}\left( T(y_{j})\right)
\right\vert \leq d_{p}^{MSL}\left( T\right) g_{p}^{L}(z),
\end{equation*}%
where $z=\sum_{j=1}^{m}\delta _{\left( x_{j},y_{j}\right) }\boxtimes
f_{j}\in Y\boxtimes W^{\#}.$ Then%
\begin{eqnarray*}
\left\vert \left\langle T,z\right\rangle \right\vert  &=&\left\vert
\sum_{j=1}^{m}f_{j}\left( T(x_{j})\right) -f_{j}\left( T(y_{j})\right)
\right\vert  \\
&\leq &d_{p}^{MSL}\left( T\right) (\sum_{j=1}^{m}\left\Vert \delta _{\left(
x_{j},y_{j}\right) }\right\Vert ^{p})^{\frac{1}{p}}\left\Vert \left(
f_{j}\right) _{j=1}^{m}\right\Vert _{\ell _{p^{\ast }}^{m,w}\left(
W^{\#}\right) } \\
&\leq &d_{p}^{MSL}\left( T\right) (\sum_{j=1}^{m}d\left( x_{j},y_{j}\right)
^{p})^{\frac{1}{p}}\left\Vert \left( f_{j}\right) _{j=1}^{m}\right\Vert
_{\ell _{p^{\ast }}^{m,w}\left( W^{\#}\right) }
\end{eqnarray*}%
\noindent $2)\Rightarrow 3):$ We will show that $\delta _{W}\circ
T:Y\rightarrow \mathcal{F}\left( W\right) $ is strongly Lipschitz $p$-summing%
\textit{. }Let\textit{\ }$\left( x_{j}\right) _{j=1}^{m},\left( y_{j}\right)
_{j=1}^{m}$\textit{\ }in $Y$ and $\left( f_{j}\right) _{j=1}^{m}\subset
W^{\#}\left( =\mathcal{F}\left( W\right) ^{\ast }\right) .$ Then%
\begin{eqnarray*}
\left\vert \sum_{j=1}^{m}\left\langle \delta _{W}\circ T(x_{j})-\delta
_{W}\circ T(y_{j}),f_{j}\right\rangle \right\vert  &=&\left\vert
\sum_{j=1}^{m}f_{j}\left( \delta _{W}\circ T(x_{j})\right) -f_{j}\left(
\delta _{W}\circ T(y_{j})\right) \right\vert  \\
&=&\left\vert \sum_{j=1}^{m}f_{j}\left( T(x_{j})\right) -f_{j}\left(
T(y_{j})\right) \right\vert  \\
&\leq &K(\sum_{j=1}^{m}d\left( x_{j},y_{j}\right) ^{p})^{\frac{1}{p}%
}\left\Vert \left( f_{j}\right) _{j=1}^{m}\right\Vert _{\ell _{p^{\ast
}}^{m,w}\left( W^{\#}\right) }
\end{eqnarray*}%
Then $\delta _{Y}\circ T$ is strongly Lipschitz $p$-summing.

\noindent $3)\Rightarrow 4):$ We know that 
\begin{equation*}
\widehat{\delta _{W}\circ T}=\widetilde{T}.
\end{equation*}%
Furthermore, according to \cite[Proposition 3.1]{13}, the linearization $%
\widehat{\delta _{W}\circ T}$ is strongly $p$-summing, which implies that $%
\widetilde{T}$ is also strongly $p$-summing.

\noindent $4)\Rightarrow 1):$ Let $z=\sum_{l=1}^{k}\delta _{\left(
x_{l},y_{l}\right) }\boxtimes f_{l}\in Y\boxtimes W^{\#}.$ Assuming that $%
\widetilde{T}$ is strongly $p$-summing, we can deduce from \cite[Proposition
6.12]{12} that $\widetilde{T}$ satisfies the following%
\begin{equation*}
\left\vert \sum_{l=1}^{k}\left\langle \widetilde{T}\left( \mathfrak{n}%
_{l}\right) ,f_{l}\right\rangle \right\vert \leq Kg_{p}(\sum_{l=1}^{k}%
\mathfrak{n}_{l}\otimes f_{l}),
\end{equation*}%
If we put $\mathfrak{n}_{l}=\delta _{\left( x_{l},y_{l}\right) }$ for $1\leq
l\leq k,$ we find%
\begin{eqnarray*}
\left\vert \sum_{l=1}^{l}\left\langle \widetilde{T}\left( \delta _{\left(
x_{l},y_{l}\right) }\right) ,f_{l}\right\rangle \right\vert  &=&\left\vert
\sum_{k=1}^{l}f_{l}\left( T\left( x_{l}\right) \right) -f_{l}\left( T\left(
y_{l}\right) \right) \right\vert  \\
&\leq &Kg_{p}(z)=Kg_{p}^{L}(z).\qquad \blacksquare 
\end{eqnarray*}

Using the same reasoning as in Corollary 2.6, we can establish the following
result.\newline

\textbf{Corollary 3.3}. \textit{Consider }$1\leq p\leq \infty $\textit{\ and
a Lipschitz operator }$T:Y\rightarrow W$\textit{\ between pointed metric
spaces. The following properties are equivalent.}

\noindent \textit{1) }$T$\textit{\ is MS-strongly Lipschitz }$p$-summing%
\textit{.}

\noindent \textit{2) The Lipschitz adjoint }$T^{\#}:W^{\#}\rightarrow Y^{\#}$
\textit{is }$p^{\ast }$\textit{-summing.}\newline

The following integral characterization is an adaptation of the linear case.
To prove it, we rely on the fact that $T^{\#}$ is $p^{\ast }$-summing or $%
\widetilde{T}$ is strongly $p$-summing.\newline

\textbf{Theorem 3.4}. \textit{Consider }$1\leq p\leq \infty $\textit{\ and a
Lipschitz operator }$T:Y\rightarrow W$\textit{\ between pointed metric
spaces. The following properties are equivalent.}

\noindent \textit{1) }$T$\textit{\ is MS-strongly Lipschitz }$p$-\textit{%
summing.}

\noindent \textit{2) There exist a constant }$K>0$ \textit{and} \textit{a} 
\textit{Radon} \textit{probability }$\mu $\textit{\ on }$B_{Lip_{0}\left(
W\right) ^{\ast }}$ \textit{such that for every }$x,y\in Y$ \textit{and} $%
f\in W^{\#},$ \textit{we have}%
\begin{equation*}
\left\vert f\left( T(x)\right) -f\left( T(y)\right) \right\vert \leq
Kd\left( x,y\right) \left( \int_{B_{Lip_{0}\left( W\right) ^{\ast
}}}\left\vert \left\langle f,\mathfrak{n}\right\rangle \right\vert ^{p^{\ast
}}d\mu \left( \mathfrak{n}\right) \right) ^{\frac{1}{p^{\ast }}}.
\end{equation*}

Let us now recall the definition of the tensor norm $w_{p}$ on the product
of two Banach spaces $F\otimes G$, which has been studied in \cite[p. 180]%
{12}. Let $p\in \left[ 1,\infty \right] \ $we have%
\begin{equation*}
w_{p}\left( z\right) =\inf \left\{ \left\Vert \left( x_{j}\right)
_{j=1}^{m}\right\Vert _{l_{p}^{m,w}\left( F\right) }\left\Vert \left(
y_{j}\right) _{j=1}^{m}\right\Vert _{l_{p^{\ast }}^{m,w}\left( G\right)
}\right\} ,
\end{equation*}%
where the infimum is taken over all representations of $z$ of the form $%
z=\sum_{j=1}^{m}x_{j}\otimes y_{j}\in F\otimes G.$\newline

\textbf{Definition 3.5}. Consider a Lipschitz operator $T:Y\rightarrow W$\
between pointed metric spaces. For $1\leq p\leq \infty ,$\textit{\ }$T$ is
said to be\ MS-Lipschitz $p$-nuclear if there is a constant $K>0$\ such that
such that for every $z=\sum_{l=1}^{k}\delta _{\left( x_{l},y_{l}\right)
}\boxtimes f_{l}\in Y\boxtimes W^{\#}$ we have 
\begin{equation}
\left\vert \left\langle T,z\right\rangle \right\vert \leq Kw_{p}^{L}(z), 
\tag{3.3}
\end{equation}%
We define $\mathcal{N}_{p}^{MSL}\left( Y,W\right) $ as the set of all
MS-Lipschitz $p$-nuclear operators from $Y$ into $W$. Additionally, $%
n_{p}^{MSL}\left( T\right) $ represents the smallest constant $K$ that
satisfies (3.3). If $W$ is a Banach space, it is important to reiterate that 
$\mathcal{N}_{p}^{MSL}\left( Y,W\right) $ does not possess the structure of
a vector space.\newline

Suppose that $W$ is a Banach space. By restricting the previous definition
to the linear forms of $W^{\ast }$, we obtain the definition of (Cohen)
Lipschitz $p$-nuclear operators as introduced in \cite{9}. Indeed, let $%
z=\sum_{l=1}^{k}\delta _{\left( x_{l},y_{l}\right) }\boxtimes a_{l}^{\ast
}\in Y\boxtimes W^{\ast }$ we have 
\begin{eqnarray}
\left\vert \left\langle T,z\right\rangle \right\vert  &=&\left\vert
\sum_{l=1}^{k}a_{l}^{\ast }\left( T\left( x_{l}\right) \right) -a_{l}^{\ast
}\left( T\left( y_{l}\right) \right) \right\vert =\left\vert
\sum_{l=1}^{k}\left\langle T\left( x_{l}\right) -T\left( y_{l}\right)
,a_{l}^{\ast }\right\rangle \right\vert   \notag \\
&\leq &Kw_{p}^{L}(z)=Kw_{p}(\sum_{l=1}^{k}\delta _{\left( x_{l},y_{l}\right)
}\otimes a_{l}^{\ast }) \\
&\leq &K\left\Vert \left( \delta _{\left( x_{l},y_{l}\right) }\right)
_{l=1}^{k}\right\Vert _{\ell _{p}^{k,w}\left( \mathcal{F}\left( Y\right)
\right) }\left\Vert \left( a_{l}^{\ast }\right) _{l=1}^{k}\right\Vert _{\ell
_{p^{\ast }}^{k,w}\left( W^{\ast }\right) } \\
&\leq &K\sup_{f\in Y^{\#}}(\sum_{l=1}^{k}\left\vert
f(x_{l})-f(y_{l})\right\vert ^{p})^{\frac{1}{p}}\sup_{\left\Vert
a\right\Vert _{W}=1}(\sum_{l=1}^{k}\left\vert \left\langle a_{l}^{\ast
},a\right\rangle \right\vert ^{p^{\ast }})^{\frac{1}{p^{\ast }}},
\end{eqnarray}

Consequently, it follows that every MS-Lipschitz $p$-nuclear operator is
also Lipschitz $p$-nuclear.\newline

\textbf{Theorem 3.6}. \textit{Consider }$1\leq p\leq \infty $\textit{\ and a
Lipschitz operator }$T:Y\rightarrow W$\textit{\ between pointed metric
spaces. The following properties are equivalent.}

\noindent \textit{1) }$T$\textit{\ is MS-Lipschitz }$p$\textit{-nuclear.}

\noindent \textit{2) There is a constant }$K>0$\textit{\ such that for every 
}$\left( x_{j}^{i}\right) _{j=1}^{m_{1}},\left( y_{j}^{i}\right)
_{j=1}^{m_{1}}$\textit{\ in }$Y,$ $\left( f_{j}\right) _{j=1}^{m_{1}}\subset
W^{\#};\left( 1\leq i\leq m_{2}\right) $ \textit{and} $m_{1},m_{2}\in 
\mathbb{N}^{\ast }$\textit{, we have}%
\begin{equation}
\left\vert \sum_{j=1}^{m_{1}}\sum_{i=1}^{m_{2}}f_{j}\left(
T(x_{j}^{i})\right) -f_{j}\left( T(y_{j}^{i})\right) \right\vert \leq
K\sup_{s\in Y^{\#}}(\sum_{j=1}^{m_{1}}\left\vert
\sum_{i=1}^{m_{2}}s(x_{j}^{i})-s(y_{j}^{i})\right\vert ^{p})^{\frac{1}{p}%
}\left\Vert \left( f_{j}\right) _{j=1}^{m_{1}}\right\Vert _{\ell _{p^{\ast
}}^{m_{1},w}\left( W^{\#}\right) }.  \tag{3.4}
\end{equation}

\noindent \textit{3) }$\widetilde{T}:\mathcal{F}\left( Y\right) \rightarrow 
\mathcal{F}\left( W\right) $\textit{\ is }$p$\textit{-nuclear.}\newline

\textbf{Proof}.

\noindent $1)\Rightarrow 2):$ Let $T$ be a MS-Lipschitz $p$-nuclear
operator. Let $\left( x_{j}^{i}\right) _{j=1}^{m_{1}},\left(
y_{j}^{i}\right) _{j=1}^{m_{1}}$\textit{\ }in $Y$ and $\left( f_{j}\right)
_{j=1}^{m_{1}}$ in $W^{\#};\left( 1\leq i\leq m_{2}\right) $\textit{, }we
have%
\begin{equation*}
\left\vert \sum_{j=1}^{m_{1}}\sum_{i=1}^{m_{2}}f_{j}\left(
T(x_{j}^{i})\right) -f_{j}\left( T(y_{j}^{i})\right) \right\vert \leq
Kg_{p}^{L}(z),
\end{equation*}%
where $z=\sum_{j=1}^{m_{1}}\sum_{i=1}^{m_{2}}\delta _{\left(
x_{j}^{i},y_{j}^{i}\right) }\boxtimes f_{j}\in Y\boxtimes W^{\#}.$ Then%
\begin{eqnarray*}
&&\left\vert \sum_{j=1}^{m_{1}}\sum_{i=1}^{m_{2}}f_{j}\left(
T(x_{j}^{i})\right) -f_{j}\left( T(y_{j}^{i})\right) \right\vert  \\
&\leq &n_{p}^{MSL}\left( T\right) \sup_{s\in
Y^{\#}}(\sum_{j=1}^{m_{1}}\left\vert s(\sum_{i=1}^{m_{2}}\delta _{\left(
x_{j}^{i},y_{j}^{i}\right) })\right\vert ^{p})^{\frac{1}{p}}\left\Vert
\left( f_{j}\right) _{j=1}^{m_{1}}\right\Vert _{\ell _{p^{\ast
}}^{m_{1},w}\left( W^{\#}\right) } \\
&\leq &n_{p}^{MSL}\left( T\right) \sup_{s\in
Y^{\#}}(\sum_{j=1}^{m_{1}}\left\vert
\sum_{i=1}^{m_{2}}s(x_{j}^{i})-s(y_{j}^{i})\right\vert ^{p})^{\frac{1}{p}%
}\left\Vert \left( f_{j}\right) _{j=1}^{m_{1}}\right\Vert _{\ell _{p^{\ast
}}^{m_{1},w}\left( W^{\#}\right) }.
\end{eqnarray*}%
\noindent $2)\Rightarrow 3):$ Let\textit{\ }$\left( \mathfrak{n}_{j}\right)
_{j=1}^{m_{1}}\subset \mathcal{F}\left( Y\right) $ and $\left( f_{j}\right)
_{j=1}^{m_{1}}\subset W^{\#}$ such that 
\begin{equation*}
\mathfrak{m}_{j}=\sum_{i=1}^{m_{2}}\delta _{\left(
x_{j}^{i},y_{j}^{i}\right) }\in \mathcal{F}\left( Y\right) ,\text{ }\left(
1\leq j\leq m_{1}\right) .
\end{equation*}%
Then%
\begin{eqnarray*}
(\sum_{j=1}^{m_{1}}\left\vert \left\langle \widetilde{T}\left( \mathfrak{n}%
_{j}\right) ,f_{j}\right\rangle \right\vert ^{p})^{\frac{1}{p}}
&=&(\sum_{j=1}^{m_{1}}\left\vert \left\langle \sum_{i=1}^{m_{2}}\left(
\delta _{W}\circ T\left( x_{j}^{i}\right) -\delta _{W}\circ T\left(
x_{j}^{i}\right) \right) ,f_{j}\right\rangle \right\vert ^{p})^{\frac{1}{p}}
\\
&=&\sum_{j=1}^{m_{1}}(\left\vert \sum_{i=1}^{m_{2}}f_{j}\left(
T(x_{j}^{i})\right) -f_{j}\left( T(y_{j}^{i})\right) \right\vert ^{p})^{%
\frac{1}{p}} \\
&\leq &K\sup_{s\in Y^{\#}}(\sum_{j=1}^{m_{1}}\left\vert
\sum_{i=1}^{m_{2}}s(x_{j}^{i})-s(y_{j}^{i})\right\vert ^{p})^{\frac{1}{p}%
}\left\Vert \left( f_{j}\right) _{j=1}^{m_{1}}\right\Vert _{\ell _{p^{\ast
}}^{m_{1},w}\left( W^{\#}\right) } \\
&\leq &K\sup_{s\in Y^{\#}}(\sum_{j=1}^{m_{1}}\left\vert
s(\sum_{i=1}^{m_{2}}\delta _{(x_{j}^{i},y_{j}^{i})})\right\vert ^{p})^{\frac{%
1}{p}}\left\Vert \left( f_{j}\right) _{j=1}^{m_{1}}\right\Vert _{\ell
_{p^{\ast }}^{m_{1},w}\left( W^{\#}\right) } \\
&\leq &K\sup_{s\in Y^{\#}}(\sum_{j=1}^{m_{1}}\left\vert s\left( \mathfrak{n}%
_{j}\right) \right\vert ^{p})^{\frac{1}{p}}\left\Vert \left( f_{j}\right)
_{j=1}^{m_{1}}\right\Vert _{\ell _{p^{\ast }}^{m_{1},w}\left( W^{\#}\right) }
\end{eqnarray*}%
Then $\widetilde{T}$\textit{\ }is $p$-nuclear.

\noindent $3)\Rightarrow 1):$ Now, we suppose that $\widetilde{T}$\textit{\ }%
is $p$-nuclear. Let $x_{l},y_{l}\in Y$ and $f_{l}\in W^{\#}$ $\left( 1\leq
l\leq k\right) $ we have 
\begin{eqnarray*}
\left\vert \sum_{l=1}^{k}f_{l}\left( T(x_{l})\right) -f_{l}\left(
T(y_{l})\right) \right\vert &=&\left\vert \sum_{l=1}^{k}\left\langle \delta
_{W}\circ T\left( x_{l}\right) -\delta _{W}\circ T\left( y_{l}\right)
,f_{l}\right\rangle \right\vert \\
&=&\left\vert \sum_{l=1}^{k}\left\langle \widehat{T}\left( \delta _{\left(
x_{l},y_{l}\right) }\right) ,f_{l}\right\rangle \right\vert \\
&\leq &n_{p}\left( \widetilde{T}\right) w_{p}(z)=n_{p}\left( \widetilde{T}%
\right) w_{p}^{L}(z),
\end{eqnarray*}%
where $z=\sum_{l=1}^{k}\delta _{\left( x_{l},y_{l}\right) }\boxtimes
f_{l}\in Y\boxtimes W^{\#}.$ Finally, $T$ is MS-Lipschitz $p$-nuclear and we
have%
\begin{equation*}
n_{p}^{MSL}\left( T\right) \leq n_{p}\left( \widetilde{T}\right) .\qquad
\blacksquare
\end{equation*}

By utilizing the result presented in \cite[Theorem 2.2.4]{3}, we can
establish the following relationship between $T$ and its Lipschitz adjoint $%
T^{\#}$.\newline

\textbf{Corollary 3.7}. \textit{Consider }$1\leq p\leq \infty $\textit{\ and
a Lipschitz operator }$T:Y\rightarrow W$\textit{\ between pointed metric
spaces. The following properties are equivalent.}

\noindent \textit{1) }$T$\textit{\ is MS-Lipschitz }$p$\textit{-nuclear.}

\noindent \textit{2) The Lipschitz adjoint }$T^{\#}:W^{\#}\rightarrow Y^{\#}$
\textit{is }$p^{\ast }$\textit{-nuclear.}\newline

The following integral characterization is an adaptation of the linear case.
Its proof will be omitted.\newline

\textbf{Theorem 3.8}. \textit{Consider }$1\leq p\leq \infty $\textit{\ and a
Lipschitz operator }$T:Y\rightarrow W$\textit{\ between pointed metric
spaces. The following properties are equivalent.}

\noindent \textit{1) }$T$\textit{\ is MS-Lipschitz }$p$\textit{-nuclear.}

\noindent \textit{2) There exist a constant }$K>0,$ \textit{a} \textit{Radon}
\textit{probability }$\mu $\textit{\ on }$B_{Y^{\#}}$ \textit{and} $\eta \in
B_{Lip_{0}\left( W\right) ^{\ast }}$ \textit{such that for every }$\left(
x^{i}\right) _{i=1}^{m},\left( y^{i}\right) _{i=1}^{m}$\textit{\ in }$Y$ and 
$f\in W^{\#},$ \textit{we have}%
\begin{eqnarray*}
\left\vert \sum_{i=1}^{m}f\left( T(x^{i})\right) -f\left( T(x^{i})\right)
\right\vert  &\leq &K(\int_{B_{Y^{\#}}}\left\vert
\sum_{i=1}^{m}s(x^{i})-s(y^{i})\right\vert ^{p}d\mu \left( s\right) )^{\frac{%
1}{p}}\times  \\
&&(\int_{B_{Lip_{0}\left( W\right) ^{\ast }}}\left\vert \left\langle s,%
\mathfrak{n}\right\rangle \right\vert ^{p^{\ast }}d\eta \left( \mathfrak{n}%
\right) )^{\frac{1}{p^{\ast }}}.
\end{eqnarray*}

\textbf{Theorem 3.9}. \textit{Let }$Y,W$\textit{\ and }$Z$\textit{\ be three
pointed metric spaces. Let }$u:Y\rightarrow W$\textit{\ be an MS-Lipschitz }$%
p$\textit{-summing operator and }$v:W\rightarrow Z$\textit{\ be an
MS-strongly Lipschitz }$p$\textit{-summing operator. Then the composition }$%
T=v\circ u$\textit{\ is MS-Lipschitz }$p$\textit{-nuclear.}\newline

\textbf{Proof}. By \cite[p. 124]{7}, we have%
\begin{equation*}
\widetilde{T}=\widetilde{v}\circ \widetilde{u}.
\end{equation*}%
According to a result due to Cohen \cite{3}, the linear operator $\widetilde{%
u}$ being $p$-summing and $\widetilde{v}$ being strongly $p$-summing imply
that $\widetilde{T}$ is $p$-nuclear. Consequently, $T$ is also MS-Lipschitz $%
p$-nuclear.\qquad $\blacksquare $\newline

\textbf{Remark 3.10}.

In the linear case, the converse of the previous statement is true. However,
in our case, it is unknown whether every MS-Lipschitz $p$-nuclear operator
can be expressed as the product of an MS-Lipschitz $p$-summing operator and
an MS-strongly $p$-summing operator.\newline

\textbf{Declarations}\newline

\textbf{Conflict of interest.} The authors declare that they have no
conflicts of interest.\newline

\end{document}